\documentclass[reqno]{amsart}

\newtheorem{Theorem}{Theorem}[section]

\newtheorem{Corollary}[Theorem]{Corollary}
\theoremstyle{remark}
\newtheorem{Remark}[Theorem]{Remark}
\numberwithin{equation}{section}

\allowdisplaybreaks

\begin{document}

\title[Inversion of bilateral series]
{Inversion of bilateral basic hypergeometric series}
\author{Michael Schlosser}

\address{Institut f\"ur Mathematik der Universit\"at Wien,
Strudlhofgasse 4, A-1090 Wien, Austria}
\email{schlosse@ap.univie.ac.at}
\urladdr{http://www.mat.univie.ac.at/{\textasciitilde}schlosse}
\thanks{The author was supported by an APART grant of the Austrian
Academy of Sciences}
\date{June 4, 2002}

\dedicatory{{\rm\small Institut f\"ur Mathematik der Universit\"at Wien,\\
Strudlhofgasse 4, A-1090 Wien, Austria\\
E-mail: \tt schlosse@ap.univie.ac.at\\
WWW: \tt http://www.mat.univie.ac.at/{\textasciitilde}schlosse}}

\begin{abstract}
We present a new matrix inverse with applications in the
theory of bilateral basic hypergeometric series. Our matrix
inversion result is directly extracted from an instance of
Bailey's very-well-poised ${}_6\psi_6$ summation theorem, and
involves two infinite matrices which are not lower-triangular.
We combine our bilateral matrix inverse with known basic
hypergeometric summation theorems to derive, via inverse relations,
several new identities for bilateral basic hypergeometric series.
\end{abstract}

\maketitle

\section{Introduction}\label{sec0}

Bailey's~\cite[Eq.~(4.7)]{bail66} very-well-poised ${}_6\psi_6$
summation formula,
\begin{multline}\label{66gl}
{}_6\psi_6\!\left[\begin{matrix}q\sqrt{a},-q\sqrt{a},b,c,d,e\\
\sqrt{a},-\sqrt{a},aq/b,aq/c,aq/d,aq/e\end{matrix}\,;q,
\frac{a^2q}{bcde}\right]\\
=\frac {(q,aq,q/a,aq/bc,aq/bd,aq/be,aq/cd,aq/ce,aq/de;q)_{\infty}}
{(aq/b,aq/c,aq/d,aq/e,q/b,q/c,q/d,q/e,a^2q/bcde;q)_{\infty}},
\end{multline}
where $|a^2q/bcde|<1$ (cf.\ \cite[Eq.~(5.3.1)]{grhyp}),
stands on the top of the classical hierarchy of summation theorems for
bilateral basic hypergeometric series. It contains many important
identities as special cases, among them Jacobi's triple product
identity, the $q$-Pfaff--Saalsch\"utz summation, and the
$q$-binomial theorem, to name just a few.
Various applications of Bailey's ${}_6\psi_6$ summation exist in
number theory (see Andrews~\cite[pp.~461--468]{andappl}) and in
special functions (see, e.g., Ismail and Masson~\cite{ismasson}).
A combinatorial (partition theoretic) application of
Bailey's ${}_6\psi_6$ summation formula was recently revealed in
remarkable work of Alladi, Andrews, and Berkovich~\cite{fourpar}.

Different proofs of \eqref{66gl} are known. A very elegant proof
using analytic continuation was given by Askey and Ismail~\cite{askmail}.
For an elementary proof using manipulations of series,
see Schlosser~\cite{schlelsum}.

In addition to Bailey's ${}_6\psi_6$ summation formula, there is a
significant number of other important summation and transformation
theorems for basic hypergeometric series (cf.\ \cite{grhyp}).
Basic hypergeometric series (and, more generally, $q$-series)
have various applications in combinatorics, number theory,
representation theory, statistics, and physics,
see Andrews~\cite{andappl},\cite{qandrews}.  For a general account
of the importance of basic hypergeometric series in the theory of
special functions see Andrews, Askey, and Roy~\cite{sfaar}.

Various techniques have been employed for the study of
basic hypergeometric series. A fundamental approach is to
start with simple identities and build up the theory by
successively deriving more complicated identities.
The heart piece of this method is the ``Bailey transform'',
a simple but efficient interchange of summation argument.
The Bailey transform is even more powerful if it is combined with
a specific summation theorem, in which case it becomes a
``Bailey lemma'' (see Andrews~\cite{qandrews},\cite{andbt}).
Starting with an identity, the Bailey lemma generates
an infinite chain (or lattice) of identities,
a so-called ``Bailey chain'', or more general, a ``Bailey lattice''
(see \cite{qandrews} and \cite{aab}).

Another important tool for proving or deriving identities is
using ``inverse relations'', which are an immediate consequence
of matrix inversions. By this method, the proof of a given identity
may be reduced to the proof of another ``dual'' identity.
On the other hand, given a known identity, by applying 
inverse relations a possibly new identity may be derived.

Matrix inversion and Bailey lemma are not unrelated tools for deriving
identities. Already Andrews~\cite{andrews} had observed that the specific
application of the Bailey transform in the classical Bailey lemma is
equivalent to an explicit matrix inversion result. Important matrix inversions
have been found by Gould and Hsu~\cite{gouldhsu}, Carlitz~\cite{carlitz},
Gessel and Stanton~\cite{gesstant}, Bressoud~\cite{bress}, Al-Salam and
Verma~\cite{alsv}, Gasper~\cite{gas}, Krattenthaler~\cite{kratt-mi},
and Warnaar~\cite{warnell}. Similar results in higher dimensions
(related to multiple series) have been obtained by Chu~\cite{chu},
Milne~\cite{milne}, Lilly and Milne~\cite{lilmil}, Bhatnagar and
Milne~\cite{bhatmil}, Schlosser~\cite{schlossmmi},\cite{schlossnmmi},
and Krattenthaler and Schlosser~\cite{krattschl}.

In this paper, we provide yet another explicit pair of infinite matrices
being inverses of each other. The difference from all the previously
mentioned matrix inverses is that our new result involves two infinite
matrices which are not necessarily lower-triangular, i.e., all their
entries may be non-zero. The corresponding orthogonality relations
are infinite convergent sums. Our new ``bilateral matrix inverse''
in Theorem~\ref{bmi} is directly extracted from an instance of
Bailey's very-well-poised ${}_6\psi_6$ summation \eqref{66gl}
and extends Bressoud's~\cite{bress} matrix inverse (which involves
lower-triangular matrices) by an additional free parameter.

After a short introduction to basic hypergeometric series and inverse
relations in Section~\ref{secpre}, we state and prove our main result,
a ``bilateral matrix inverse'', in Section~\ref{secbmi}. In
Section~\ref{secappl}, we combine our new matix inverse with
basic hypergeometric summation theorems to derive, via inverse relations,
new bilateral summation theorems. Finally, in Section~\ref{secmore},
we use our newly derived summations from Section~\ref{secappl} to deduce
by elementary manipulations of sums further bilateral series identities.

In a forthcoming paper, we apply part of the current analysis to multiple
sums. In particular, by appropriately specializing Gustafson's~\cite{gus}
$A_r$ and $C_r$ ${}_6\psi_6$ summations, we derive multidimensional
extensions of our bilateral matrix inverse in Theorem~\ref{bmi},
and deduce some multilateral summations as applications.

We are currently preparing another article which features a new
bilateral Bailey lemma, based on our new bilateral matrix inverse in
Theorem~\ref{bmi}, combined with Bailey's~\cite[Eq.~(3.3.)]{bailwp}
nonterminating $_8\phi_7$ summation. This bilateral Bailey lemma is
different from (and does not specialize to) the Bailey lemmas
considered in \cite{qandrews} or \cite{andbt}.

This paper has been typeset while the author was visiting Northwestern
University in Evanston, Illinois, for the Spring Quarter 2002.
The author wishes to acknowledge the positive research atmosphere
experienced there. We are in particular thankful for stimulating
discussions with George Gasper.

\section{Preliminaries}\label{secpre}
\subsection{Notation and basic hypergeometric series}

Here we recall some standard notation for $q$-series,
and basic hypergeometric series (cf.\ \cite{grhyp}).

Let $q$ be a complex number such that $0<|q|<1$. We define the
{\em $q$-shifted factorial} for all integers $k$ by
\begin{equation*}
(a;q)_{\infty}:=\prod_{j=0}^{\infty}(1-aq^j)\qquad\text{and}\qquad
(a;q)_k:=\frac{(a;q)_{\infty}}{(aq^k;q)_{\infty}}.
\end{equation*}
For brevity, we employ the condensed notation
\begin{equation*}
(a_1,\ldots,a_m;q)_k\equiv (a_1;q)_k\dots(a_m;q)_k
\end{equation*}
where $k$ is an integer or infinity. Further, we utilize
\begin{equation}\label{defhyp}
{}_s\phi_{s-1}\!\left[\begin{matrix}a_1,a_2,\dots,a_s\\
b_1,b_2,\dots,b_{s-1}\end{matrix}\,;q,z\right]:=
\sum _{k=0} ^{\infty}\frac {(a_1,a_2,\dots,a_s;q)_k}
{(q,b_1,\dots,b_{s-1};q)_k}z^k,
\end{equation}
and
\begin{equation}\label{defhypb}
{}_s\psi_s\!\left[\begin{matrix}a_1,a_2,\dots,a_s\\
b_1,b_2,\dots,b_s\end{matrix}\,;q,z\right]:=
\sum _{k=-\infty} ^{\infty}\frac {(a_1,a_2,\dots,a_s;q)_k}
{(b_1,b_2,\dots,b_s;q)_k}z^k,
\end{equation}
to denote the {\em basic hypergeometric ${}_s\phi_{s-1}$ series},
and the {\em bilateral basic hypergeometric ${}_s\psi_s$ series},
respectively. In \eqref{defhyp} or \eqref{defhypb}, $a_1,\dots,a_s$ are
called the {\em upper parameters}, $b_1,\dots,b_s$ the
{\em lower parameters}, $z$ is the {\em argument}, and 
$q$ the {\em base} of the series.
See \cite[p.~25 and p.~125]{grhyp} for the criteria
of when these series terminate, or, if not, when they converge. 
 
The classical theory of basic hypergeometric series contains
numerous summation and transformation formulae
involving ${}_s\phi_{s-1}$ or ${}_s\psi_s$ series.
Many of these summation theorems require
that the parameters satisfy the condition of being
either balanced and/or very-well-poised.
An ${}_s\phi_{s-1}$ basic hypergeometric series is called
{\em balanced} if $b_1\cdots b_{s-1}=a_1\cdots a_sq$ and $z=q$.
An ${}_s\phi_{s-1}$ series is {\em well-poised} if
$a_1q=a_2b_1=\cdots=a_sb_{s-1}$. An ${}_s\phi_{s-1}$ basic
hypergeometric series is called {\em very-well-poised}
if it is well-poised and if $a_2=-a_3=q\sqrt{a_1}$.
Note that the factor
\begin{equation*}
\frac {1-a_1q^{2k}}{1-a_1}
\end{equation*}
appears in  a very-well-poised series.
The parameter $a_1$ is usually referred to as the
{\em special parameter} of such a series.
Similarly, a bilateral ${}_s\psi_s$ basic hypergeometric series is
well-poised if $a_1b_1=a_2b_2\cdots=a_sb_s$ and very-well-poised if,
in addition, $a_1=-a_2=qb_1=-qb_2$. Further, we call a bilateral
${}_s\psi_s$ basic hypergeometric series balanced if
$b_1\cdots b_s=a_1\cdots a_sq^2$ and $z=q$.

A standard reference for basic hypergeometric series
is Gasper and Rahman's text~\cite{grhyp}.
In our computations in the subsequent sections
we frequently use some elementary identities of
$q$-shifted factorials, listed in \cite[Appendix~I]{grhyp}.

In the following we display some summation theorems which we utilize
in Sections~\ref{secappl} and \ref{secmore}.

One of the most important theorems in the theory of basic
hypergeometric series is 
Jackson's~\cite{jacksum} terminating very-well-poised balanced
${}_8\phi_7$ summation (cf.\ \cite[Eq.~(2.6.2)]{grhyp}):
\begin{multline}\label{87gl}
{}_8\phi_7\!\left[\begin{matrix}a,\,q\sqrt{a},-q\sqrt{a},b,c,d,
a^2q^{1+n}/bcd,q^{-n}\\
\sqrt{a},-\sqrt{a},aq/b,aq/c,aq/d,bcdq^{-n}/a,aq^{1+n}\end{matrix}\,;q,
q\right]\\
=\frac {(aq,aq/bc,aq/bd,aq/cd;q)_n}
{(aq/b,aq/c,aq/d,aq/bcd;q)_n}.
\end{multline}

A less well known but nevertheless very useful identity is
the following very-well-poised ${}_8\phi_7$ summation:
\begin{multline}\label{qwatsongl}
{}_8\phi_7\!\left[\begin{matrix}\lambda,\,q\sqrt{\lambda},-q\sqrt{\lambda},
a,b,c,-c,\lambda q/c^2\\
\sqrt{\lambda},-\sqrt{\lambda},\lambda q/a,\lambda q/b,\lambda q/c,
-\lambda q/c,c^2\end{matrix}\,;q,
-\frac{\lambda q}{ab}\right]\\
=\frac {(\lambda q,c^2/\lambda;q)_{\infty}
(aq,bq,c^2q/a,c^2q/b;q^2)_{\infty}}
{(\lambda q/a,\lambda q/b;q)_{\infty}
(q,abq,c^2q,c^2q/ab;q^2)_{\infty}},
\end{multline}
provided $|\lambda q/ab|<1$, where $\lambda =-c\sqrt{ab/q}$.

The ${}_8\phi_7$ summation formula in \eqref{qwatsongl} is a $q$-analogue
of a ${}_3F_2$ summation due to Whipple~\cite{whipple1} (but
commonly attributed to Watson~\cite{watson} who gave that ${}_3F_2$
summation for the case when the series is terminating).

A bilateral summation even (slightly) more general than the ${}_6\psi_6$
sum in \eqref{66gl} is H.~S.~Shukla's~\cite{shukla} very-well-poised
${}_8\psi_8$ summation:
\begin{multline}\label{shuklagl}
{}_8\psi_8\!\left[\begin{matrix}q\sqrt{a},-q\sqrt{a},b,c,d,e,f,aq^2/f\\
\sqrt{a},-\sqrt{a},aq/b,aq/c,aq/d,aq/e,aq/f,f/q\end{matrix}\,;q,
\frac{a^2}{bcde}\right]\\
=\left(1-\frac{(1-bc/a)(1-bd/a)(1-be/a)}
{(1-bq/f)(1-bf/aq)(1-bcde/a^2)}\right)
\frac{(1-f/bq)(1-bf/aq)}{(1-f/aq)(1-f/q)}\\\times
\frac {(q,aq,q/a,aq/bc,aq/bd,aq/be,aq/cd,aq/ce,aq/de;q)_{\infty}}
{(aq/b,aq/c,aq/d,aq/e,q/b,q/c,q/d,q/e,a^2q/bcde;q)_{\infty}},
\end{multline}
where $|a^2/bcde|<1$.
Note that \eqref{shuklagl} reduces to \eqref{66gl} if $f\to 0$
or $f\to\infty$. For a generalization of \eqref{shuklagl},
see \eqref{milnekm2gl}.

We will use the summations \eqref{87gl}, \eqref{qwatsongl}, and
\eqref{shuklagl} in Section~\ref{secappl}
to derive new bilateral summation theorems.

\subsection{Inverse relations}

Let $\mathbb Z$ denote the set of integers.
In the following, we consider infinite matrices
$(f_{nk})_{n,k\in\mathbb Z}$ and
$(g_{nk})_{n,k\in\mathbb Z}$, and infinite sequences
$(a_n)_{n\in\mathbb Z}$ and $(b_n)_{n\in\mathbb Z}$.

We say that the infinite matrices 
$(f_{nk})_{n,k\in\mathbb Z}$ and $(g_{kl})_{k,l\in\mathbb Z}$ are
{\it inverses} of each other if and only if the following orthogonality
relation holds:
\begin{equation}\label{orthrel}
\sum _{k\in\mathbb Z}f_{nk}g_{kl}=\delta_{nl}\qquad\qquad
\text {for all}\quad n,l\in\mathbb Z.
\end{equation}
Clearly, since inverse matrices commute, we also then have
\begin{equation}\label{orthreld}
\sum _{l\in\mathbb Z}g_{kl}f_{lj}=\delta_{kj}\qquad\qquad
\text {for all}\quad k,j\in\mathbb Z.
\end{equation}

Note that in \eqref{orthrel} and  \eqref{orthreld} we are {\em not}
requiring that the infinite matrices are lower-triangular.
If they were, the summations on the left hand sides of \eqref{orthrel} and
\eqref{orthreld} would be in fact finite sums. In the general case,
the sums will be infinite. If the summands of the infinite series involve
complex numbers, we require suitable convergence conditions to hold
(such as absolute convergence; for interchanging
double sums we also need uniform convergence).

It is immediate from the orthogonality relations \eqref{orthrel}
and  \eqref{orthreld} that the following {\em inverse relations} hold:
Let $(f_{nk})_{n,k\in\mathbb Z}$ and $(g_{kl})_{k,l\in\mathbb Z}$
be infinite matrices being inverses of each other. Then
\begin{equation}\label{invrel1}
\sum _{k\in\mathbb Z}f_{n k}a_{k}=b_{n}\qquad\qquad\text {for all $n$,}
\end{equation}
if and only if
\begin{equation}\label{invrel2}
\sum _{l\in\mathbb Z}g_{k l}b_{l}=a_{k}\qquad\qquad\text {for all $k$.}
\end{equation}

The other variant of inverse relations, which may be called
``rotated inversion" (see also Riordan~\cite{riord}), reads as follows:
Let $(f_{nk})_{n,k\in\mathbb Z}$ and $(g_{kl})_{k,l\in\mathbb Z}$
be infinite matrices being inverses of each other. Then
\begin{equation}\label{rotinv1}
\sum _{n\in\mathbb Z}f_{n k}a_{n}=b_{k}\qquad\qquad\text {for all $k$,}
\end{equation}
if and only if
\begin{equation}\label{rotinv2}
\sum _{k\in\mathbb Z}g_{k l}b_{k}=a_{l}\qquad\qquad\text {for all $l$.}
\end{equation}

We also note here that if the considered sequences
involve complex numbers we need suitable convergence conditions for
the above inverse relations \eqref{invrel1}/\eqref{invrel2}
and \eqref{rotinv1}/\eqref{rotinv2} to hold. 

Inverse relations are a powerful tool for proving or deriving
identities. For instance, given an identity in the form \eqref{rotinv2},
we can immediately deduce \eqref{rotinv1}, which may possibly be a new
identity. It is exactly this variant of inverse relations which we will
utilize in Section~\ref{secappl} to derive new summation theorems for
bilateral series.

\section{A new bilateral matrix inverse}
\label{secbmi}

We now present our main result, an explicit pair of inverse infinite
matrices which are not lower-triangular.

\begin{Theorem}\label{bmi}
Let $a$, $b$, and $c$ be indeterminates. The infinite matrices
$(f_{nk})_{n,k\in\mathbb Z}$ and $(g_{kl})_{k,l\in\mathbb Z}$ are {\em
inverses} of each other where
\begin{multline}\label{bmf}
f_{nk}=\frac{(aq/b,bq/a,aq/c,cq/a,bq,q/b,cq,q/c;q)_{\infty}}
{(q,q,aq,q/a,aq/bc,bcq/a,cq/b,bq/c;q)_{\infty}}\\
\times\frac{(1-bcq^{2n}/a)}{(1-bc/a)}\frac{(b;q)_{n+k}\,(a/c;q)_{k-n}}
{(cq;q)_{n+k}\,(aq/b;q)_{k-n}}
\end{multline}
and
\begin{equation}\label{bmg}
g_{kl}=\frac{(1-aq^{2k})}{(1-a)}\frac{(c;q)_{k+l}\,(a/b;q)_{k-l}}
{(bq;q)_{k+l}\,(aq/c;q)_{k-l}}\,q^{k-l}.
\end{equation}
\end{Theorem}

\begin{Remark}
If we let $c\to a$ in Theorem~\ref{bmi}, we obtain a matrix inverse
found by Bressoud~\cite{bress} which he directly extracted from the
terminating very-well-poised ${}_6\phi_5$ summation
(a special case of \eqref{66gl}).
If, after letting $c\to a$, we additionally let $a\to 0$, we obtain
Andrews'~\cite[Lemma~3]{andrews} ``Bailey transform matrices'', a
matrix inversion underlying the powerful Bailey lemma.
\end{Remark}

\begin{proof}[Proof of Theorem~\ref{bmi}]
We show that the inverse matrices \eqref{bmf}/\eqref{bmg} satisfy
the orthogonality relation \eqref{orthrel}.
Writing out the sum $\sum _{k\in\mathbb Z}f_{nk}g_{kl}$ with the above
choices of $f_{nk}$ and $g_{kl}$ we observe that
the series can be summed by an application of
Bailey's very-well-poised ${}_6\psi_6$ summation \eqref{66gl}.
The specializations needed there are $b\mapsto bq^n$, $c\mapsto cq^l$,
$d\mapsto aq^{-l}/b$, and $e\mapsto aq^{-n}/c$. Bailey's formula then
gives us a product containing the factors
$(q^{1-n+l},q^{1+n-l};q)_{\infty}$. Since
\begin{equation*}
(q^{1-n+l},q^{1+n-l};q)_{\infty}=0
\end{equation*}
for all integers $n$ and $l$ with $n\neq l$, we can simplify
the product (setting $n=l$, the only non-zero case)
and readily determine that the sum indeed boils down to $\delta_{nl}$.
The details are as follows:
\begin{multline*}
\sum _{k\in\mathbb Z}f_{nk}g_{kl}
=\sum _{k=-\infty}^{\infty}
\frac{(aq/b,bq/a,aq/c,cq/a,bq,q/b,cq,q/c;q)_{\infty}}
{(q,q,aq,q/a,aq/bc,bcq/a,cq/b,bq/c;q)_{\infty}}\\\times
\frac{(1-bcq^{2n}/a)}{(1-bc/a)}\frac{(b;q)_{n+k}(a/c;q)_{k-n}}
{(cq;q)_{n+k}(aq/b;q)_{k-n}}
\frac{(1-aq^{2k})}{(1-a)}\frac{(c;q)_{k+l}(a/b;q)_{k-l}}
{(bq;q)_{k+l}(aq/c;q)_{k-l}}q^{k-l}\\
=\frac{(aq/b,bq/a,aq/c,cq/a,bq,q/b,cq,q/c;q)_{\infty}}
{(q,q,aq,q/a,aq/bc,bcq/a,cq/b,bq/c;q)_{\infty}}\\\times
\frac{(1-bcq^{2n}/a)}{(1-bc/a)}\frac{(b;q)_n(a/c;q)_{-n}(c;q)_l(a/b;q)_{-l}}
{(cq;q)_n(aq/b;q)_{-n}(bq;q)_l(aq/c;q)_{-l}}q^{-l}\\\times
\sum _{k=-\infty}^{\infty}\frac{(1-aq^{2k})}{(1-a)}
\frac{(bq^n,aq^{-n}/c,cq^l,aq^{-l}/b;q)_k}
{(aq^{1-n}/b,cq^{1+n},aq^{1-l}/c,bq^{1+l};q)_k}q^k\\
=\frac{(aq/b,bq/a,aq/c,cq/a,bq,q/b,cq,q/c;q)_{\infty}}
{(q,q,aq,q/a,aq/bc,bcq/a,cq/b,bq/c;q)_{\infty}}\\\times
\frac{(1-bcq^{2n}/a)}{(1-bc/a)}\frac{(b;q)_n(a/c;q)_{-n}(c;q)_l(a/b;q)_{-l}}
{(cq;q)_n(aq/b;q)_{-n}(bq;q)_l(aq/c;q)_{-l}}q^{-l}\\\times
\frac {(q,aq,q/a,cq/b,aq^{1-n-l}/bc,q^{1-n+l},q^{1+n-l},
bcq^{1+n+l}/a,bq/c;q)_{\infty}}
{(aq^{1-n}/b,cq^{1+n},aq^{1-l}/c,bq^{1+l},
q^{1-n}/b,cq^{1+n}/a,q^{1-l}/c,bq^{1+l}/a,q;q)_{\infty}}\\
=\delta_{nl},
\end{multline*}
where we have been using some elementary identities
involving $q$-shifted factorials (cf.\ \cite[Appendix~I]{grhyp}).
\end{proof}

Observe that the dual orthogonality relation \eqref{orthreld}
of the matrices \eqref{bmf}/\eqref{bmg} is established
by a similar instance of Bailey's very-well-poised $_6\psi_6$ summation
(where the argument is again $q$). In particular, both series
in \eqref{orthrel} and in \eqref{orthreld} converge absolutely
for $|q|<1$.

\begin{Remark}
We do not need the full ${}_6\psi_6$ summation \eqref{66gl}
to prove the above orthogonality relation(s). All we really need
is the very-well-poised {\em balanced} ${}_6\psi_6$ summation,
\begin{multline}\label{66agl}
\sum _{k=-\infty}^{\infty}\frac{(1-aq^{2k})}{(1-a)}
\frac{(b,c,d,a^2/bcd;q)_k}
{(aq/b,aq/c,aq/d,bcdq/a;q)_k}q^k\\
=\frac {(aq,q/a,aq/bc,bcq/a,aq/bd,bdq/a,cdq/a,aq/cd;q)_{\infty}}
{(aq/b,aq/c,aq/d,q/b,q/c,q/d,bcdq/a,bcdq/a^2;q)_{\infty}},
\end{multline}
which is just a special case of Bailey's very-well-poised $_6\psi_6$
summation \eqref{66gl}. Usually, behind a particular orthogonality
relation there is a more general summation theorem. In case
of finite sums, this summation theorem can often be proved by a simple
telescoping argument. A similar situation happens here with \eqref{66agl}.
Gasper and Rahman~\cite{gasrahmind} showed that the indefinite bibasic sum
\begin{multline*}
\sum_{k=-m}^n\frac{(1-adp^kq^k)(1-bp^k/dq^k)}{(1-ad)(1-b/d)}
\frac{(a,b;p)_k(c,ad^2/bc;q)_k}{(dq,adq/b;q)_k(adp/c,bcp/d;p)_k}q^k\\
=\frac{(1-a)(1-b)(1-c)(1-ad^2/bc)}{d(1-ad)(1-b/d)(1-c/d)(1-ad/bc)}\\\times
\left(\frac{(ap,bp;p)_n(cq,ad^2q/bc;q)_n}{(dq,adq/b;q)_n(adp/c,bcp/d;p)_n}
-\frac{(c/ad,d/bc;p)_{m+1}(1/d,b/ad;q)_{m+1}}
{(1/c,bc/ad^2;q)_{m+1}(1/a,1/b;p)_{m+1}}\right)
\end{multline*}
telescopes. Replacing $a$ and $d$, respectively, by $d$ and $a/d$,
letting $n,m\to\infty$, and setting $p=q$, the above summation reduces to
\begin{multline}\label{66bgl}
\sum _{k=-\infty}^{\infty}\frac{(1-aq^{2k})}{(1-a)}
\frac{(b,c,d,a^2/bcd;q)_k}
{(aq/b,aq/c,aq/d,bcdq/a;q)_k}q^k\\
=\frac{(1-b)(1-c)(1-d)(1-bcd/a^2)}{(1-a)(1-bd/a)(1-cd/a)(1-bc/a)}\\\times
\left(\frac{(bq,cq,dq,a^2q/bcd;q)_{\infty}}
{(aq/b,aq/c,aq/d,bcdq/a;q)_{\infty}}
-\frac{(b/a,c/a,d/a,a/bcd;q)_{\infty}}
{(1/b,1/c,1/d,bcd/a^2;q)_{\infty}}\right).
\end{multline}
Now, (miraculously) the right hand side of \eqref{66bgl} can be
transformed into the right hand side of \eqref{66agl} using a theta
function identity of Weierstrass~\cite[p.~451, Example 5]{whitwatson}
(also cf.\ \cite[Ex.~5.21]{grhyp}).
\end{Remark}

\section{Some applications}\label{secappl}

Here we combine the matrix inverse of Theorem~\ref{bmi}
with specific summation theorems to derive, via inverse relations,
new identites for bilateral basic hypergeometric series.
For convenience, we only use the (rotated) inverse relations
\eqref{rotinv1}/\eqref{rotinv2}, but because of the symmetric
structure of our inverse matrices we could as well also
employ the inverse relations \eqref{invrel1}/\eqref{invrel2}
to obtain equivalent results.

In our first application of Theorem~\ref{bmi} we apply inverse relations
to Jackson's terminating very-well-poised balanced ${}_8\phi_7$ summation,
and obtain a summation for a particular very-well-poised balanced
${}_8\psi_8$ series, see Theorem~\ref{88n}.
In our second application we apply inverse relations
to the nonterminating ${}_8\phi_7$ summation in \eqref{qwatsongl}.
This leads us to a new bilateral quadratic summation, see Theorem~\ref{qun}.
Whereas in our first two applications of Theorem~\ref{bmi} we invert a
terminating summation, and a nonterminating unilateral summation,
respectively, our third application of Theorem~\ref{bmi} involves inverting
a genuine bilateral summation theorem.
It turns out that if we invert the ${}_6\psi_6$ summation theorem
\eqref{66gl}, we again end up with \eqref{66gl}. (We may also invert
Rogers'~\cite[p.~29, second eq.]{rogers} nonterminating ${}_6\phi_5$
summation (cf.\ \cite[Eq.~(2.7.1)]{grhyp}) which leads to a summation
for a ${}_6\psi_6$ series. Unfortunately the ${}_6\psi_6$ summation
obtained in this way is not as general as \eqref{66gl}.)
In order to derive something new, we climb up higher in the hierarchy of
bilateral basic hypergeometric series and invert Shukla's very-well-poised
${}_8\psi_8$ summation \eqref{shuklagl}. The  result is Theorem~\ref{shn}.

We start with our first application. 
Using some elementary identities for $q$-shifted factorials on the
right hand side, we may deduce from Jackson's terminating very-well-poised
balanced ${}_8\phi_7$ summation \eqref{87gl} the identity
\begin{multline}\label{87gla}
{}_8\phi_7\!\left[\begin{matrix}a,\,q\sqrt{a},-q\sqrt{a},d,cq^l,
aq^{-l}/b,abq^{1+N}/cd,q^{-N}\\
\sqrt{a},-\sqrt{a},aq/d,aq^{1-l}/c,bq^{1+l},cdq^{-N}/b,aq^{1+N}\end{matrix}\,
;q,q\right]\\
=\frac {(aq,aq/cd,bq/c,bq/d;q)_N}
{(bq/cd,bq,aq/d,aq/c;q)_N}
\frac{(bq,bq^{1+N}/d,cd/a,cq^{-N}/a;q)_l}
{(c/a,cdq^{-N}/a,bq/d,bq^{1+N};q)_l},
\end{multline}
where $N$ is a nonnegative integer.
Note that we have introduced an additional integer $l$ in the above
${}_8\phi_7$ summation. We can rewrite \eqref{87gla} as
\begin{multline}\label{ex1gl}
\sum_{k\in\mathbb Z}\underset{g_{kl}}
{\underbrace{\frac{(1-aq^{2k})}{(1-a)}\frac{(c;q)_{k+l}(a/b;q)_{k-l}}
{(bq;q)_{k+l}(aq/c;q)_{k-l}}q^{k-l}}}\;
\underset{b_k}{\underbrace{\frac{(a,d,abq^{1+N}/cd,q^{-N};q)_k}
{(q,aq/d,cdq^{-N}/b,aq^{1+N};q)_k}}}\\
=\underset{a_l}{\underbrace{
\frac {(aq,aq/cd,bq/c,bq/d;q)_N}
{(bq/cd,bq,aq/d,aq/c;q)_N}
\frac{(bq^{1+N}/d,c,cd/a,cq^{-N}/a;q)_l}
{(cdq^{-N}/a,bq/a,bq/d,bq^{1+N};q)_l}\left(\frac bc\right)^l}},
\end{multline}
by which we have established \eqref{rotinv2} for the above choices
of $a_l$ and $b_k$. By virtue of the matrix inversion in Theorem~\ref{bmi}
and the equivalence of \eqref{rotinv1} and \eqref{rotinv2}, we
immediately deduce from \eqref{ex1gl} the following inverse relation:
\begin{multline*}
\sum_{n\in\mathbb Z}\underset{f_{nk}}
{\underbrace{\frac{(aq/b,bq/a,aq/c,cq/a,bq,q/b,cq,q/c;q)_{\infty}
(1-bcq^{2n}/a)(b;q)_{n+k}(a/c;q)_{k-n}}
{(q,q,aq,q/a,aq/bc,bcq/a,cq/b,bq/c;q)_{\infty}
(1-bc/a)(cq;q)_{n+k}(aq/b;q)_{k-n}}}}\\\times
\underset{a_n}{\underbrace{\frac {(aq,aq/cd,bq/c,bq/d;q)_N}
{(bq/cd,bq,aq/d,aq/c;q)_N}
\frac {(bq^{1+N}/d,c,cd/a,cq^{-N}/a;q)_n}{(cdq^{-N}/a,bq/a,bq/d,bq^{1+N};q)_n}
\left(\frac bc\right)^n}}\\
=\underset{b_k}{\underbrace{\frac{(a,d,abq^{1+N}/cd,q^{-N};q)_k}
{(q,aq/d,cdq^{-N}/b,aq^{1+N};q)_k}}}.
\end{multline*}
Now, after substitution of variables
(simultaneously $a\mapsto cd/a$, $b\mapsto d$, and $d\mapsto bd/a$),
and moving some factors to the other side, we obtain
the following very-well-poised balanced ${}_8\psi_8$ summation:

\begin{Theorem}\label{88n}
Let $a$, $b$, $c$, and $d$ be indeterminates,
let $k$ be an arbitrary integer and $N$ a nonnegative integer. Then
\begin{multline}\label{88ngl}
{}_8\psi_8\!\left[\begin{matrix}q\sqrt{a},-q\sqrt{a},b,c,dq^k,
aq^{-k}/c,aq^{1+N}/b,aq^{-N}/d\\
\sqrt{a},-\sqrt{a},aq/b,aq/c,aq^{1-k}/d,
cq^{1+k},bq^{-N},dq^{1+N}\end{matrix}\,;q,q\right]\\
=\frac {(aq/bc,cq/b,dq,dq/a;q)_N}
{(cdq/a,dq/c,q/b,aq/b;q)_N}
\frac{(cd/a,bd/a,cq,cq/a,dq^{1+N}/b,q^{-N};q)_k}
{(q,cq/b,d/a,d,bcq^{-N}/a,cdq^{1+N}/a;q)_k}\\\times
\frac{(q,q,aq,q/a,cdq/a,aq/cd,cq/d,dq/c;q)_{\infty}}
{(cq,q/c,dq,q/d,cq/a,aq/c,dq/a,aq/d;q)_{\infty}}.
\end{multline}
\end{Theorem}
Note that two of the upper parameters of the ${}_8\psi_8$ series
in \eqref{88ngl} differ multiplicatively from corresponding lower
parameters by $q^N$, a nonnegative integral power of $q$.

Although, to the best of our knowledge, Theorem~\ref{88n}
is stated here explicitly for the first time, there is also
another way to derive (or verify) \eqref{88ngl}
which we just sketch very briefly.
By adequately specializing
M.~Jackson's~\cite[Eq.~(2.2)]{mjack1} transformation
formula for a very-well-poised $_8\psi_8$ series into a sum of two
$_8\phi_7$ series (cf.\ \cite[Eq.~(5.6.2)]{grhyp}),
one of the right hand side terms
(a multiple of an $_8\phi_7$ series) becomes zero,
while the other $_8\phi_7$ series can be summed
by an application of F.~H.~Jackson's terminating $_8\phi_7$
summation. We omit displaying the details.

Our second application involves the nonterminating
very-well-poised ${}_8\phi_7$ summation in \eqref{qwatsongl}
which, by introducing an additional integer $l$, can be written as
\begin{multline}\label{qwatsagl}
{}_8\phi_7\!\left[\begin{matrix}a,\,q\sqrt{a},-q\sqrt{a},cq^l,
aq^{-l}/b,-(abq/c)^{1/2},(abq/c)^{1/2},c/b\\
\sqrt{a},-\sqrt{a},aq^{1-l}/c,bq^{1+l},
-(acq/b)^{1/2},(acq/b)^{1/2},abq/c\end{matrix}\,;q,
-\frac{bq}c\right]\\
=\frac{(aq,bq/c;q)_{\infty}}{(aq^{1-l}/c,bq^{1+l};q)_{\infty}}
\frac{(cq^{1+l},aq^{1-l}/b,abq^{2-l}/c^2,b^2q^{2+l}/c;q^2)_{\infty}}
{(q,acq/b,abq^2/c,b^2q^2/c^2;q^2)_{\infty}},
\end{multline}
where $|bq/c|<1$. Note that on the right hand side
we have $q$ and $q^2$ appearing as bases. Thus, we may refer to
\eqref{qwatsagl} as a {\em quadratic} summation.
We may rewrite \eqref{qwatsagl} as
\begin{multline}\label{ex2gl}
\sum_{k\in\mathbb Z}\underset{g_{kl}}
{\underbrace{\frac{(1-aq^{2k})}{(1-a)}\frac{(c;q)_{k+l}(a/b;q)_{k-l}}
{(bq;q)_{k+l}(aq/c;q)_{k-l}}q^{k-l}}}\;
\underset{b_k}{\underbrace{\frac{(a,c/b;q)_k}{(q,abq/c;q)_k}\,
\frac{(abq/c;q^2)_k}{(acq/b;q^2)_k}\left(-\frac bc\right)^k}}\\
=\underset{a_l}{\underbrace{\frac{(c;q)_l(a/b;q)_{-l}
(aq,bq/c;q)_{\infty}
(cq^{1+l},aq^{1-l}/b,abq^{2-l}/c^2,
b^2q^{2+l}/c;q^2)_{\infty}}
{q^l\,(aq/c,bq;q)_{\infty}
(q,acq/b,abq^2/c,b^2q^2/c^2;q^2)_{\infty}}}},
\end{multline}
by which we have established \eqref{rotinv2} for the above choices
of $a_l$ and $b_k$. By virtue of the matrix inversion in Theorem~\ref{bmi}
and the equivalence of \eqref{rotinv1} and \eqref{rotinv2}, we immediately
deduce from \eqref{ex2gl} the following inverse relation:
\begin{multline*}
\sum_{n\in\mathbb Z}\underset{f_{nk}}
{\underbrace{
\frac{(aq/b,bq/a,aq/c,cq/a,bq,q/b,cq,q/c;q)_{\infty}
(1-bcq^{2n}/a)(b;q)_{n+k}(a/c;q)_{k-n}}
{(q,q,aq,q/a,aq/bc,bcq/a,cq/b,bq/c;q)_{\infty}
(1-bc/a)(cq;q)_{n+k}(aq/b;q)_{k-n}}}}\\\times
\underset{a_n}
{\underbrace{\frac{(c;q)_n(a/b;q)_{-n}
(aq,bq/c;q)_{\infty}(cq^{1+n},aq^{1-n}/b,abq^{2-n}/c^2,
b^2q^{2+n}/c;q^2)_{\infty}}
{q^n\,(aq/c,bq;q)_{\infty}
(q,acq/b,abq^2/c,b^2q^2/c^2;q^2)_{\infty}}}}\\
=\underset{b_k}
{\underbrace{\frac{(a,c/b;q)_k}{(q,abq/c;q)_k}
\frac{(abq/c;q^2)_k}{(acq/b;q^2)_k}\left(-\frac bc\right)^k}}
\end{multline*}
(subject to the convergence condition $|q|<1$.)
After substitution of variables (simultaneously $a\mapsto bc/a$,
$b\mapsto c$, and $c\mapsto b$), we obtain the following
bilateral quadratic summation:

\begin{Theorem}\label{qun}
Let $a$, $b$, and $c$ be indeterminates,
and let $k$ be an arbitrary integer. Then
\begin{multline}\label{qungl}
\sum_{n=-\infty}^{\infty}
\frac{(1-aq^{2n})}{(1-a)}\frac{(b,cq^k,aq^{-k}/b;q)_n}
{(aq/b,aq^{1-k}/c,bq^{1+k};q)_n}\,
q^{\binom n2}\left(-\frac{aq}c\right)^n\\\times
(bq^{1+n},bq^{1-n}/a,c^2q^{2+n}/b,
c^2q^{2-n}/ab;q^2)_{\infty}\\
=\frac{(bc/a,bq,bq/a,b/c;q)_k}
{(q,c/a,c,c^2q/a;q)_k}\,
\frac{(c^2q/a;q^2)_k}{(b^2q/a;q^2)_k}
\left(-\frac cb\right)^k\\\times
\frac{(q,q,aq,q/a,aq/bc,bq/c;q)_{\infty}}
{(bq,q/b,bq/a,aq/b,q/c,aq/c;q)_{\infty}}\,
(q,b^2q/a,c^2q^2/a,c^2q^2/b^2;q^2)_{\infty}.
\end{multline}
\end{Theorem}
Note that the series in \eqref{qungl} converges absolutely for $|q|<1$.
This can easily be seen by splitting the sum in two series
according to the parity of the index of summation ($n=2m$ or
$n=2m+1$ for integer $m$). In this case it is also easy to
rewrite the whole sum as a sum of two hypergeometric series
with base $q^2$, where each sum is absolutly convergent for
$|q^2|<1$ (i.e., $|q|<1$).

Our third application involves Shukla's very-well-poised
${}_8\psi_8$ summation \eqref{shuklagl}. By introducing an additional
integer $l$, we can deduce from \eqref{shuklagl} the identity
\begin{multline}\label{ex3gl}
{}_8\psi_8\!\left[\begin{matrix}q\sqrt{a},-q\sqrt{a},
cq^l,aq^{-l}/b,d,e,u,aq^2/u\\
\sqrt{a},-\sqrt{a},aq^{1-l}/c,bq^{1+l},aq/d,
aq/e,aq/u,u/q\end{matrix}\,;q,
\frac{ab}{cde}\right]\\
=\left(1-\frac{(1-c/b)(1-cdq^l/a)(1-ceq^l/a)}
{(1-cq^{1+l}/u)(1-cuq^{l-1}/a)(1-cde/ab)}\right)\\\times
\frac{(1-uq^{-1-l}/c)(1-cuq^{l-1}/a)}
{(1-u/aq)(1-u/q)}\left(\frac a{de}\right)^l
\frac{(cd/a,ce/a,bq,bq/a;q)_l}{(c/a,c,bq/d,bq/e;q)_l}\\\times
\frac {(q,aq,q/a,bq/c,aq/cd,aq/ce,bq/d,bq/e,aq/de;q)_{\infty}}
{(aq/c,bq,aq/d,aq/e,q/c,bq/a,q/d,q/e,abq/cde;q)_{\infty}},
\end{multline}
where $|ab/cde|<1$. We may rewrite \eqref{ex3gl} in the form \eqref{rotinv2}
with $g_{kl}$ as in \eqref{bmg},
\begin{multline*}
a_l=\left(1-\frac{(1-c/b)(1-cdq^l/a)(1-ceq^l/a)}
{(1-cq^{1+l}/u)(1-cuq^{l-1}/a)(1-cde/ab)}\right)\\\times
\frac{(1-uq^{-1-l}/c)(1-cuq^{l-1}/a)}
{(1-u/aq)(1-u/q)}\left(\frac{ab}{cde}\right)^l
\frac{(cd/a,ce/a;q)_l}{(bq/d,bq/e;q)_l}\\\times
\frac {(q,aq,q/a,bq/c,aq/cd,aq/ce,bq/d,bq/e,aq/de;q)_{\infty}}
{(aq/c,bq,aq/d,aq/e,q/c,bq/a,q/d,q/e,abq/cde;q)_{\infty}},
\end{multline*}
and
\begin{equation*}
b_k=\frac{(d,e,u,aq^2/u;q)_k}{(aq/d,aq/e,aq/u,u/q;q)_k}
\left(\frac{ab}{cdeq}\right)^k.
\end{equation*}
This implies the inverse relation \eqref{rotinv1} with $f_{nk}$ as in
\eqref{bmf}, and the above choices of $a_n$ and $b_k$.
This inverse relation is equivalent to
\begin{multline}\label{ex4gl}
\sum_{n=-\infty}^{\infty}\frac{(1-bcq^{2n}/a)}{(1-bc/a)}
\frac{(bq^k,bq^{-k}/a,cdq/a,ceq/a;q)_n}
{(cq^{1-k}/a,cq^{1+k},bq/d,bq/e;q)_n}
\left(\frac a{de}\right)^n\\\times
\left(1-\frac{(1-cde/ab)(1-cq^{n+1}/u)(1-cuq^{n-1}/a)}
{(1-c/b)(1-cdq^n/a)(1-ceq^n/a)}\right)\\
=\frac{(1-uq^{k-1})(1-aq^{k+1}/u)}{(b-c)}
\left(\frac{ab}{cdeq}\right)^k
\frac{(cq,aq/b,d,e;q)_k}
{(a/c,b,aq/d,aq/e;q)_k}\\\times
\frac {(q,aq/d,aq/e,q/d,q/e,cq/b,aq/bc,bcq/a,ab/cde;q)_{\infty}}
{(a/cd,a/ce,bq/d,bq/e,q/b,aq/b,cq,cq/a,aq/de;q)_{\infty}},
\end{multline}
provided $|a/de|<1$. Now, it is not difficult to see that the integer
parameter $k$ in \eqref{ex4gl} can be eliminated. Hence, we could
readily set $k=0$. However, our following
substitutions not only serve to eliminate the superfluous $k$
but also serve to obtain more symmetry of the parameters.
Specifically, in \eqref{ex4gl} we perform the simultaneous
substitutions $a\mapsto bq^{-2k}/c$, $b\mapsto bq^{-k}$,
$c\mapsto aq^{-k}/c$, $d\mapsto bdq^{-k}/a$, $e\mapsto beq^{-k}/a$,
$u\mapsto aq^{1-k}/cu$, simplify further and obtain:

\begin{Theorem}\label{shn}
Let $a$, $b$, $c$, $d$, $e$, and $u$ be indeterminates. Then
\begin{multline}\label{shngl}
\sum_{n=-\infty}^{\infty}\frac{(1-aq^{2n})}{(1-a)}
\frac{(b,c;q)_n}{(aq/b,aq/c;q)_n}
\frac{(dq,eq;q)_n}{(aq/d,aq/e;q)_n}
\left(\frac{a^2}{bcde}\right)^n\\\times
\left(1-\frac{(1-de/a)(1-uq^n)(1-a^2q^n/bcu)}
{(1-a/bc)(1-dq^n)(1-eq^n)}\right)\\
=\frac{(1-a/cu)(1-bu/a)
(q,aq,q/a,aq/bc,aq/bd,aq/be,aq/cd,aq/ce,a/de;q)_{\infty}}
{(b-a/c)(q/b,aq/b,q/c,aq/c,1/d,aq/d,1/e,aq/e,a^2q/bcde;q)_{\infty}},
\end{multline}
provided $|a^2/bcde|<1$.
\end{Theorem}
Note that Theorem~\ref{shn} reduces to Bailey's very-well-poised
${}_6\psi_6$ summation \eqref{66gl} if we let $u\to e$, and then
perform the substitution $e\mapsto e/q$.

\begin{Remark}
Instead of inverting Shukla's very-well-poised ${}_8\psi_8$ summation
\eqref{shuklagl}, we could also have inverted a more general identity
to obtain an even more general result. A good candidate is the $k=2$
special case of a ${}_{2k+4}\psi_{2k+4}$ transformation formula
due to Milne~\cite[Theorem~1.7]{milnekm} which can be written as follows:
\begin{multline}\label{milnekm2gl}
{}_8\psi_8\!\left[\begin{matrix}q\sqrt{a},-q\sqrt{a},b,c,d,e,f,aq^{1+N}/f\\
\sqrt{a},-\sqrt{a},aq/b,aq/c,aq/d,aq/e,aq/f,fq^{-N}\end{matrix}\,;q,
\frac{a^2q^{1-N}}{bcde}\right]\\
=\frac {(q,aq,q/a,aq/bc,aq/bd,aq/be,aq/cd,aq/ce,aq/de;q)_{\infty}}
{(aq/b,aq/c,aq/d,aq/e,q/b,q/c,q/d,q/e,a^2q/bcde;q)_{\infty}}\\\times
\frac{(bq/f,aq/f;q)_N}{(aq/f,q/f;q)_N}\,
{}_4\phi_3\!\left[\begin{matrix}q^{-N},bc/a,bd/a,be/a\\
bq/f,bfq^{-N}/a,bcde/a^2\end{matrix}\,;q,q\right],
\end{multline}
where $|a^2q^{1-N}/bcde|<1$, for convergence.
Note that for $N=0$ \eqref{milnekm2gl} reduces to \eqref{66gl}, 
whereas \eqref{shuklagl} is just the special case $N=1$ of \eqref{milnekm2gl}
where the balanced ${}_4\phi_3$ series on the right hand side reduces
to a sum of two terms only. Another noteworthy special case is obtained
by letting $e\to a/b$ in \eqref{milnekm2gl}. In this case, the
${}_4\phi_3$ series on the right hand side is just one and the above
transformation reduces to a summation formula for an ${}_8\psi_8$ series
(which in an equivalent form is displayed in \eqref{ex7gl}).
A similar specialization
in the full Theorem~1.7 of Milne~\cite{milnekm}, without the
restriction $k=2$, yields a summation for a $q$-IPD type (or
$q$-Karlsson--Minton type; ``IPD'' stands for {\em integral parameter
differences}, cf.\ \cite{schlkm},\cite[Sec.~7]{schleltf})
very-well-poised ${}_{2k+4}\psi_{2k+4}$ series,
a summation apparently originally missed by Milne, but later
independently (re-)discovered by Chu~\cite[Theorem~12]{chuwp}.

However, it is clear that it is always possible to obtain more general
results. The reason that in our last application of this section
we only inverted \eqref{shuklagl} and not the more general
\eqref{milnekm2gl} is that we feel that Theorem~\ref{shn},
which already  seems to be new, looks appealing enough,
in spite, or because, of the relative simplicity of the
explicitly given terms. We leave it as an easy exercise
for the reader to carry out the analysis of inverting \eqref{milnekm2gl}
to obtain a result which is more general than Theorem~\ref{shn}.
\end{Remark}

\section{More identities}\label{secmore}

We can use the newly derived summations in Section~\ref{secappl}
to derive other identities.
The machinery we apply is very elementary. In our summations
involving an integer parameter $k$, we multiply both sides
by a suitable expression depending on $k$ and then sum both
sides of the identity over all integers $k$.
We will then have a double sum on one side of the identity.
There we interchange summations and simplify the inner sum.
As result we obtain a transformation of series.

It should be understood that the identities we are deriving
by the described method (by combining various summations) are
definitely not an exhaustive list of all possibilities.
We just give a number of samples which seem to be interesting
enough to present.

\subsection{Some consequences of Theorem~\ref{88n}}

First, we combine Theorem~\ref{88n} with itself.
We multiply both sides of \eqref{88ngl} by
\begin{equation*}
\frac{(1-cdq^{2k}/a)}{(1-cd/a)}\frac{(d/a,d,e,c,cdq^{1+M}/ae,cq^{-M}/a;q)_k}
{(cq,cq/a,cdq/ae,dq/a,eq^{-M},dq^{1+M};q)_k}q^k
\end{equation*}
and sum over all integers $k$. On the right hand side we obtain
\begin{multline}\label{ex10gl}
\frac {(aq/bc,cq/b,dq,dq/a;q)_N}
{(cdq/a,dq/c,q/b,aq/b;q)_N}
\frac{(q,q,aq,q/a,cdq/a,aq/cd,cq/d,dq/c;q)_{\infty}}
{(cq,q/c,dq,q/d,cq/a,aq/c,dq/a,aq/d;q)_{\infty}}\\\times
{}_{10}\phi_9\!\left[\begin{matrix}cd/a,\,
q(cd/a)^{\frac 12},-q(cd/a)^{\frac 12},bd/a,c,e,\\
(cd/a)^{\frac 12},-(cd/a)^{\frac 12},cq/b,dq/a,cdq/ae,\end{matrix}\right.\\*
\left.\begin{matrix}cdq^{1+M}/ae,cq^{-M}/a,dq^{1+N}/b,q^{-N}\\
eq^{-M},dq^{1+M},bcq^{-N}/a,cdq^{1+N}/a\end{matrix}\,;q,
q\right].
\end{multline}
On the left hand side we obtain
\begin{multline*}
\sum_{k=-\infty}^{\infty}
\frac{(1-cdq^{2k}/a)}{(1-cd/a)}
\frac{(d/a,d,e,c,cdq^{1+M}/ae,cq^{-M}/a;q)_k}
{(cq,cq/a,cdq/ae,dq/a,eq^{-M},dq^{1+M};q)_k}q^k\\\times
\sum_{j=-\infty}^{\infty}\frac{(1-aq^{2j})}{(1-a)}
\frac{(b,c,dq^k,aq^{-k}/c,aq^{1+N}/b,aq^{-N}/d;q)_j}
{(aq/b,aq/c,aq^{1-k}/d,cq^{1+k},bq^{-N},dq^{1+N};q)_j}q^j\\
=\sum_{j=-\infty}^{\infty}\frac{(1-aq^{2j})}{(1-a)}
\frac{(b,c,d,a/c,aq^{1+N}/b,aq^{-N}/d;q)_j}
{(aq/b,aq/c,aq/d,cq,bq^{-N},dq^{1+N};q)_j}q^j\\\times
\sum_{k=-\infty}^{\infty}
\frac{(1-cdq^{2k}/a)}{(1-cd/a)}
\frac{(dq^{-j}/a,dq^j,e,c,cdq^{1+M}/ae,cq^{-M}/a;q)_k}
{(cq^{1+j},cq^{1-j}/a,cdq/ae,dq/a,eq^{-M},dq^{1+M};q)_k}q^k.
\end{multline*}
Now the inner sum can be evaluated by \eqref{88ngl} and we obtain
\begin{multline*}
\sum_{j=-\infty}^{\infty}\frac{(1-aq^{2j})}{(1-a)}
\frac{(b,c,d,a/c,aq^{1+N}/b,aq^{-N}/d;q)_j}
{(aq/b,aq/c,aq/d,cq,bq^{-N},dq^{1+N};q)_j}q^j\\\times
\frac {(dq/ae,cq/e,dq,aq/c;q)_M}
{(aq,dq/c,q/e,cdq/ae;q)_M}
\frac{(a,ae/c,cq,aq/d,dq^{1+M}/e,q^{-M};q)_j}
{(q,cq/e,a/c,d,aeq^{-M}/d,aq^{1+M};q)_j}\\\times
\frac{(q,q,cdq/a,aq/cd,aq,q/a,cq/d,dq/c;q)_{\infty}}
{(cq,q/c,dq,q/d,aq/d,dq/a,aq/c,cq/a;q)_{\infty}}\\
=\frac{(q,q,cdq/a,aq/cd,aq,q/a,cq/d,dq/c;q)_{\infty}}
{(cq,q/c,dq,q/d,aq/d,dq/a,aq/c,cq/a;q)_{\infty}}
\frac {(dq/ae,cq/e,dq,aq/c;q)_M}
{(aq,dq/c,q/e,cdq/ae;q)_M}\\\times
{}_{10}\phi_9\!\left[\begin{matrix}a,\,qa^{\frac 12},-qa^{\frac 12},
b,c,ae/c,aq^{1+N}/b,aq^{-N}/d,dq^{1+M}/e,q^{-M}\\
a^{\frac 12},-a^{\frac 12},aq/b,aq/c,cq/e,
bq^{-N},dq^{1+N},aeq^{-M}/d,aq^{1+M}\end{matrix}\,;q,q\right].
\end{multline*}
Equating the last expression with \eqref{ex10gl}, and performing the
simultaneous substitutions $d\mapsto e$ and $e\mapsto cd/a$, we obtain
the following terminating very-well-poised balanced ${}_{10}\phi_9$
transformation:
\begin{Corollary}\label{ex11}
Let $a$, $b$, $c$, $d$, and $e$, be indeterminates,
and let $N$ and $M$ be nonnegative integers. Then
\begin{multline}\label{ex11gl}
{}_{10}\phi_9\!\left[\begin{matrix}a,\,qa^{\frac 12},-qa^{\frac 12},
b,c,d,aq^{1+N}/b,aq^{-N}/e,aeq^{1+M}/cd,q^{-M}\\
a^{\frac 12},-a^{\frac 12},aq/b,aq/c,aq/d,
bq^{-N},eq^{1+N},cdq^{-M}/e,aq^{1+M}\end{matrix}\,;q,q\right]\\
=\frac{(aq,eq/c,aq/cd,eq/d;q)_M}{(eq/cd,aq/d,eq,aq/c;q)_M}
\frac {(aq/bc,cq/b,eq,eq/a;q)_N}{(ceq/a,eq/c,q/b,aq/b;q)_N}\\\times
{}_{10}\phi_9\!\left[\begin{matrix}ce/a,\,
q(ce/a)^{\frac 12},-q(ce/a)^{\frac 12},be/a,c,cd/a,\\
(ce/a)^{\frac 12},-(ce/a)^{\frac 12},cq/b,eq/a,eq/d,\end{matrix}\right.\\*
\left.\begin{matrix}eq^{1+M}/d,cq^{-M}/a,eq^{1+N}/b,q^{-N}\\
cdq^{-M}/a,eq^{1+M},bcq^{-N}/a,ceq^{1+N}/a\end{matrix}\,;q,
q\right].
\end{multline}
\end{Corollary}
It appears that Corollary~\ref{ex11} does {\em not} immediately
follow by specializing either of Bailey's~\cite{bailwp} four-term
nonterminating ${}_{10}\phi_9$ transformations
(cf.\ \cite[Eq.~(2.12.9) and Ex.~2.30]{grhyp}).

Now let us combine Theorem~\ref{88n} with the terminating Jackson
summation in \eqref{87gl}. We multiply both sides of \eqref{88ngl} by
\begin{equation*}
\frac{(1-cdq^{2k}/a)}{(1-cd/a)}
\frac{(cd/a,d/a,d,e,c^2q^{1+M}/ae,q^{-M};q)_k}
{(q,cq,cq/a,cdq/ae,deq^{-M}/c,cdq^{1+M}/a;q)_k}q^k
\end{equation*}
and sum over all integers $k$. On the right hand side we obtain
\begin{multline}\label{ex12gl}
\frac {(aq/bc,cq/b,dq,dq/a;q)_N}
{(cdq/a,dq/c,q/b,aq/b;q)_N}
\frac{(q,q,aq,q/a,cdq/a,aq/cd,cq/d,dq/c;q)_{\infty}}
{(cq,q/c,dq,q/d,cq/a,aq/c,dq/a,aq/d;q)_{\infty}}\\\times
{}_{10}\phi_9\!\left[\begin{matrix}cd/a,\,
q(cd/a)^{\frac 12},-q(cd/a)^{\frac 12},cd/a,bd/a,e,\\
(cd/a)^{\frac 12},-(cd/a)^{\frac 12},q,cq/b,cdq/ae,\end{matrix}\right.\\*
\left.\begin{matrix}dq^{1+N}/b,c^2q^{1+M}/ae,q^{-N},q^{-M}\\
bcq^{-N}/a,deq^{-M}/c,cdq^{1+N}/a,cdq^{1+M}/a\end{matrix}\,;q,
q\right].
\end{multline}
On the left hand side we obtain
\begin{multline*}
\sum_{k=-\infty}^{\infty}
\frac{(1-cdq^{2k}/a)}{(1-cd/a)}
\frac{(cd/a,d/a,d,e,c^2q^{1+M}/ae,q^{-M};q)_k}
{(q,cq,cq/a,cdq/ae,deq^{-M}/c,cdq^{1+M}/a;q)_k}q^k\\\times
\sum_{j=-\infty}^{\infty}\frac{(1-aq^{2j})}{(1-a)}
\frac{(b,c,dq^k,aq^{-k}/c,aq^{1+N}/b,aq^{-N}/d;q)_j}
{(aq/b,aq/c,aq^{1-k}/d,cq^{1+k},bq^{-N},dq^{1+N};q)_j}q^j\\
=\sum_{j=-\infty}^{\infty}\frac{(1-aq^{2j})}{(1-a)}
\frac{(b,c,d,a/c,aq^{1+N}/b,aq^{-N}/d;q)_j}
{(aq/b,aq/c,aq/d,cq,bq^{-N},dq^{1+N};q)_j}q^j\\\times
\sum_{k=-\infty}^{\infty}
\frac{(1-cdq^{2k}/a)}{(1-cd/a)}
\frac{(cd/a,dq^{-j}/a,dq^j,e,c^2q^{1+M}/ae,q^{-M};q)_k}
{(q,cq^{1+j},cq^{1-j}/a,cdq/ae,deq^{-M}/c,cdq^{1+M}/a;q)_k}q^k.
\end{multline*}
Now the inner sum can be evaluated by Jackson's
very-well-poised ${}_8\phi_7$ summation \eqref{87gl} and we obtain
\begin{multline*}
\sum_{j=-\infty}^{\infty}\frac{(1-aq^{2j})}{(1-a)}
\frac{(b,c,d,a/c,aq^{1+N}/b,aq^{-N}/d;q)_j}
{(aq/b,aq/c,aq/d,cq,bq^{-N},dq^{1+N};q)_j}q^j\\\times
\frac{(cdq/a,cq/d,cq^{1+j}/e,cq^{1-j}/ae;q)_M}
{(cq^{1+j},cq^{1-j}/a,cdq/ae,cq/de;q)_M}\\
=\frac{(cdq/a,cq/d,cq/e,cq/ae;q)_M}
{(cq,cq/a,cdq/ae,cq/de;q)_M}\\\times
{}_{10}\psi_{10}\!\left[\begin{matrix}qa^{\frac 12},-qa^{\frac 12},
b,c,d,ae/c,aq^{1+N}/b,aq^{-N}/d,cq^{1+M}/e,aq^{-M}/c\\
a^{\frac 12},-a^{\frac 12},aq/b,aq/c,aq/d,cq/e,
bq^{-N},dq^{1+N},aeq^{-M}/c,cq^{1+M}\end{matrix}\,;q,q\right].
\end{multline*}
Equating the last expression with \eqref{ex12gl}, and performing the
substitution $e\mapsto ce/a$, we obtain the following transformation
for a particular very-well-poised balanced ${}_{10}\psi_{10}$ series:
\begin{Corollary}\label{cor3}
Let $a$, $b$, $c$, $d$, and $e$, be indeterminates,
and let $N$ and $M$ be nonnegative integers. Then
\begin{multline}\label{ex13gl}
{}_{10}\psi_{10}\!\left[\begin{matrix}qa^{\frac 12},-qa^{\frac 12},
b,c,d,e,aq^{1+N}/b,aq^{-N}/d,aq^{1+M}/e,aq^{-M}/c\\
a^{\frac 12},-a^{\frac 12},aq/b,aq/c,aq/d,aq/e,
bq^{-N},dq^{1+N},eq^{-M},cq^{1+M}\end{matrix}\,;q,q\right]\\
=\frac{(cq,cq/a,dq/e,aq/de;q)_M}{(cdq/a,cq/d,aq/e,q/e;q)_M}
\frac {(aq/bc,cq/b,dq,dq/a;q)_N}{(cdq/a,dq/c,q/b,aq/b;q)_N}\\\times
{}_{10}\phi_9\!\left[\begin{matrix}cd/a,\,
q(cd/a)^{\frac 12},-q(cd/a)^{\frac 12},cd/a,bd/a,ce/a,\\
(cd/a)^{\frac 12},-(cd/a)^{\frac 12},q,cq/b,dq/e,\end{matrix}\right.\\*
\left.\begin{matrix}dq^{1+M}/b,cq^{1+M}/e,q^{-N},q^{-M}\\
bcq^{-N}/a,deq^{-M}/a,cdq^{1+N}/a,cdq^{1+M}/a\end{matrix}\,;q,
q\right].
\end{multline}
\end{Corollary}

Let us now combine Theorem~\ref{88n} with Shukla's
very-well-poised ${}_8\psi_8$ summation in \eqref{shuklagl}. We
multiply both sides of \eqref{88ngl} by
\begin{equation*}
\frac{(1-cdq^{2k}/a)}{(1-cd/a)}\frac{(d/a,d,e,f,gq,cdq/ag;q)_k}
{(cq,cq/a,cdq/ae,cdq/af,cd/ag,g;q)_k}\left(\frac{c^2}{aef}\right)^k
\end{equation*}
and sum over all integers $k$. On the right hand side we obtain
\begin{multline}\label{ex8gl}
\frac {(aq/bc,cq/b,dq,dq/a;q)_N}
{(cdq/a,dq/c,q/b,aq/b;q)_N}
\frac{(q,q,aq,q/a,cdq/a,aq/cd,cq/d,dq/c;q)_{\infty}}
{(cq,q/c,dq,q/d,cq/a,aq/c,dq/a,aq/d;q)_{\infty}}\\\times
{}_{10}\phi_9\!\left[\begin{matrix}cd/a,\,
q(cd/a)^{\frac 12},-q(cd/a)^{\frac 12},bd/a,e,f,\\
(cd/a)^{\frac 12},-(cd/a)^{\frac 12},
cq/b,cdq/ae,cdq/af,\end{matrix}\right.\\*
\left.\begin{matrix}gq,cdq/ag,dq^{1+N}/b,q^{-N}\\
cd/ag,g,bcq^{-N}/a,cdq^{1+N}/a\end{matrix}\,;q,
\frac{c^2}{aef}\right].
\end{multline}
On the left hand side we obtain
\begin{multline*}
\sum_{k=-\infty}^{\infty}
\frac{(1-cdq^{2k}/a)}{(1-cd/a)}\frac{(d/a,d,e,f,gq,cdq/ag;q)_k}
{(cq,cq/a,cdq/ae,cdq/af,cd/ag,g;q)_k}\left(\frac{c^2}{aef}\right)^k\\\times
\sum_{j=-\infty}^{\infty}\frac{(1-aq^{2j})}{(1-a)}
\frac{(b,c,dq^k,aq^{-k}/c,aq^{1+N}/b,aq^{-N}/d;q)_j}
{(aq/b,aq/c,aq^{1-k}/d,cq^{1+k},bq^{-N},dq^{1+N};q)_j}q^j\\
=\sum_{j=-\infty}^{\infty}\frac{(1-aq^{2j})}{(1-a)}
\frac{(b,c,d,a/c,aq^{1+N}/b,aq^{-N}/d;q)_j}
{(aq/b,aq/c,aq/d,cq,bq^{-N},dq^{1+N};q)_j}q^j\\\times
\sum_{k=-\infty}^{\infty}
\frac{(1-cdq^{2k}/a)}{(1-cd/a)}\frac{(dq^{-j}/a,dq^j,e,f,gq,cdq/ag;q)_k}
{(cq^{1+j},cq^{1-j}/a,cdq/ae,cdq/af,cd/ag,g;q)_k}
\left(\frac{c^2}{aef}\right)^k.
\end{multline*}
Now the inner sum, provided $|c^2/aef|<1$, can be evaluated by Shukla's
very-well-poised ${}_8\psi_8$ summation \eqref{shuklagl} and we obtain
\begin{multline*}
\sum_{j=-\infty}^{\infty}\frac{(1-aq^{2j})}{(1-a)}
\frac{(b,c,d,a/c,aq^{1+N}/b,aq^{-N}/d;q)_j}
{(aq/b,aq/c,aq/d,cq,bq^{-N},dq^{1+N};q)_j}q^j\\\times
\left(1-\frac{(1-d/c)(1-eq^{-j}/c)(1-fq^{-j}/c)}
{(1-dq^{-j}/ag)(1-gq^{-j}/c)(1-aef/c^2)}\right)
\frac{(1-agq^j/d)(1-gq^{-j}/c)}{(1-ag/cd)(1-g)}\\\times
\frac{(q,cdq/a,aq/cd,cq/d,cq^{1+j}/e,cq^{1+j}/f,
cq^{1-j}/ae,cq^{1-j}/af,cdq/aef;q)_{\infty}}
{(cq^{1+j},cq^{1-j}/a,cdq/ae,cdq/af,aq^{1+j}/d,q^{1-j}/d,q/e,q/f,
c^2q/aef;q)_{\infty}}\\
=\frac{(q,cdq/a,aq/cd,cq/d,cq/e,cq/f,cq/ae,cq/af,cdq/aef;q)_{\infty}}
{(cq,cq/a,cdq/ae,cdq/af,aq/d,q/d,q/e,q/f,c^2q/aef;q)_{\infty}}\\\times
\sum_{j=-\infty}^{\infty}\frac{(1-aq^{2j})}{(1-a)}
\frac{(b,c,ae/c,af/c,aq^{1+N}/b,aq^{-N}/d;q)_j}
{(aq/b,aq/c,cq/e,cq/f,bq^{-N},dq^{1+N};q)_j}
\left(\frac{cdq}{aef}\right)^j\\\times
\left(1-\frac{(1-d/c)(1-eq^{-j}/c)(1-fq^{-j}/c)}
{(1-dq^{-j}/ag)(1-gq^{-j}/c)(1-aef/c^2)}\right)
\frac{(1-agq^j/d)(1-gq^{-j}/c)}{(1-ag/cd)(1-g)}.
\end{multline*}
Equating the last expression with \eqref{ex8gl}, moving the infinite
products, and performing the simultaneous substitutions $d\mapsto f$,
$e\mapsto cd/a$, and $f\mapsto ce/a$, we obtain
\begin{Corollary}\label{cor2}
Let $a$, $b$, $c$, $d$, $e$, $f$, and $g$ be indeterminates,
and let $N$ be a nonnegative integer. Then
\begin{multline}\label{ex9gl}
\sum_{j=-\infty}^{\infty}\frac{(1-aq^{2j})}{(1-a)}
\frac{(b,c,d,e,aq^{1+N}/b,aq^{-N}/f;q)_j}
{(aq/b,aq/c,aq/d,aq/e,bq^{-N},fq^{1+N};q)_j}
\left(\frac{afq}{cde}\right)^j\\\times
\left(1-\frac{(1-f/c)(1-dq^{-j}/a)(1-eq^{-j}/a)}
{(1-fq^{-j}/ag)(1-gq^{-j}/c)(1-de/a)}\right)
\frac{(1-agq^j/f)(1-gq^{-j}/c)}{(1-ag/cf)(1-g)}\\
=\frac{(q,aq,q/a,aq/cd,aq/ce,aq/de,fq/c,fq/d,fq/e;q)_{\infty}}
{(aq/c,aq/d,aq/e,q/c,q/d,q/e,fq,fq/a,afq/cde;q)_{\infty}}\\\times
\frac {(aq/bc,cq/b,fq,fq/a;q)_N}{(cfq/a,fq/c,q/b,aq/b;q)_N}\,
{}_{10}\phi_9\!\left[\begin{matrix}cf/a,\,
q(cf/a)^{\frac 12},-q(cf/a)^{\frac 12},\\
(cf/a)^{\frac 12},-(cf/a)^{\frac 12},\end{matrix}\right.\\*
\left.\begin{matrix}bf/a,cd/a,ce/a,gq,cfq/ag,fq^{1+N}/b,q^{-N}\\
cq/b,fq/d,fq/e,cf/ag,g,bcq^{-N}/a,cfq^{1+N}/a\end{matrix}\,;q,
\frac{a}{de}\right],
\end{multline}
provided $|af/cde|<1$.
\end{Corollary}
Note that we have applied analytic continuation to obtain the convergence
condition in Corollary~\ref{cor2}.

If we let $g\to 0$ or $g\to\infty$ in \eqref{ex9gl}, we obtain
\begin{multline}\label{ex6gl}
{}_8\psi_8\!\left[\begin{matrix}q\sqrt{a},-q\sqrt{a},b,c,d,e,
aq^{1+N}/b,aq^{-N}/f\\
\sqrt{a},-\sqrt{a},aq/b,aq/c,aq/d,aq/e,
bq^{-N},fq^{1+N}\end{matrix}\,;q,
\frac{afq}{cde}\right]\\
=\frac{(q,aq,q/a,aq/cd,aq/ce,aq/de,fq/c,fq/d,fq/e;q)_{\infty}}
{(aq/c,aq/d,aq/e,q/c,q/d,q/e,fq,fq/a,afq/cde;q)_{\infty}}\\\times
\frac {(aq/bc,cq/b,fq,fq/a;q)_N}{(cfq/a,fq/c,q/b,aq/b;q)_N}\\\times
{}_8\phi_7\!\left[\begin{matrix}cf/a,\,q(cf/a)^{\frac 12},-q(cf/a)^{\frac 12},
bf/a,cd/a,ce/a,fq^{1+N}/b,q^{-N}\\
(cf/a)^{\frac 12},-(cf/a)^{\frac 12},cq/b,fq/d,fq/e,
bcq^{-N}/a,cfq^{1+N}/a\end{matrix}\,;q,
\frac{aq}{de}\right],
\end{multline}
provided $|afq/cde|<1$.
Note that two of the upper parameters in the ${}_8\psi_8$ on the
left hand side of \eqref{ex6gl} differ multiplicatively from corresponding
lower parameters by $q^N$, a nonnegative integral power of $q$.
If we let $e\to a/d$ in \eqref{ex6gl}, the ${}_8\phi_7$
on the right hand side becomes balanced and can be
simplified by Jackson's ${}_8\phi_7$ summation in \eqref{87gl}.
The resulting summation,
\begin{multline}\label{ex7gl}
{}_8\psi_8\!\left[\begin{matrix}q\sqrt{a},-q\sqrt{a},b,c,d,a/d,
aq^{1+N}/b,aq^{-N}/f\\
\sqrt{a},-\sqrt{a},aq/b,aq/c,aq/d,dq,
bq^{-N},fq^{1+N}\end{matrix}\,;q,
\frac{fq}{c}\right]\\
=\frac{(q,q,aq,q/a,aq/cd,dq/c,fq/d,dfq/a;q)_{\infty}}
{(dq,q/d,aq/d,dq/a,aq/c,q/c,fq,fq/a;q)_{\infty}}
\frac {(aq/bd,dq/b,fq,fq/a;q)_N}{(dfq/a,fq/d,q/b,aq/b;q)_N},
\end{multline}
is originally due to H.~S.~Shukla~\cite[p.~266, final remark]{shukla}.
Subsequently, various far-reaching extensions of \eqref{ex7gl} have been
found, see Milne~\cite[Theorem~1.7]{milnekm}, Chu~\cite[Theorem~12]{chuwp},
Schlosser~\cite[Sec.~8]{schleltf},\cite[Sec.~4]{schlkm},
and Rosengren~\cite{rosengrf},\cite{rosengrfc}.

\subsection{Some consequences of Theorem~\ref{qun}}

First, we combine Theorem~\ref{qun} with Jackson's
very-well-poised ${}_8\phi_7$ summation in \eqref{87gl}.
We multiply both sides of \eqref{qungl} by
\begin{equation*}
\frac{(1-bcq^{2k}/a)}{(1-bc/a)}\frac{(bc/a,c,c/a,d,b^2q^{1+N}/ad,q^{-N};q)_k}
{(q,bq/a,bq,bcq/ad,cdq^{-N}/b,bcq^{1+N}/a;q)_k}q^k
\end{equation*}
and sum over all integers $k$. On the right hand side we obtain
\begin{multline}\label{ex14gl}
\frac{(q,q,aq,q/a,aq/bc,bq/c;q)_{\infty}}
{(bq,q/b,bq/a,aq/b,q/c,aq/c;q)_{\infty}}\,
(q,b^2q/a,c^2q^2/a,c^2q^2/b^2;q^2)_{\infty}\\\times
\sum_{k=-\infty}^{\infty}
\frac{(1-bcq^{2k}/a)}{(1-bc/a)}
\frac{(bc/a,bc/a,b/c,d,
b^2q^{1+N}/ad,q^{-N};q)_k}
{(q,q,c^2q/a,bcq/ad,cdq^{-N}/b,bcq^{1+N}/a;q)_k}\\\times
\frac{(c^2q/a;q^2)_k}{(b^2q/a;q^2)_k}
\left(-\frac{cq}b\right)^k.
\end{multline}
On the left hand side we obtain
\begin{multline*}
\sum_{k=-\infty}^{\infty}
\frac{(1-bcq^{2k}/a)}{(1-bc/a)}\frac{(bc/a,c,c/a,d,b^2q^{1+N}/ad,q^{-N};q)_k}
{(q,bq/a,bq,bcq/ad,cdq^{-N}/b,bcq^{1+N}/a;q)_k}q^k\\\times
\sum_{n=-\infty}^{\infty}
\frac{(1-aq^{2n})}{(1-a)}\frac{(b,cq^k,aq^{-k}/b;q)_n}
{(aq/b,aq^{1-k}/c,bq^{1+k};q)_n}
q^{\binom n2}\left(-\frac{aq}c\right)^n\\\times
(bq^{1+n},bq^{1-n}/a,c^2q^{2+n}/b,
c^2q^{2-n}/ab;q^2)_{\infty}\\
=\sum_{n=-\infty}^{\infty}
\frac{(1-aq^{2n})}{(1-a)}\frac{(b,c,a/b;q)_n}{(aq/b,aq/c,bq;q)_n}
q^{\binom n2}\left(-\frac{aq}c\right)^n\\\times
(bq^{1+n},bq^{1-n}/a,c^2q^{2+n}/b,c^2q^{2-n}/ab;q^2)_{\infty}\\\times
\sum_{k=-\infty}^{\infty}
\frac{(1-bcq^{2k}/a)}{(1-bc/a)}\frac{(bc/a,cq^n,cq^{-n}/a,d,
b^2q^{1+N}/ad,q^{-N};q)_k}
{(q,bq^{1-n}/a,bq^{1+n},bcq/ad,cdq^{-N}/b,bcq^{1+N}/a;q)_k}q^k.
\end{multline*}
Now the inner sum can be evaluated by \eqref{87gl} and we obtain
\begin{multline*}
\sum_{n=-\infty}^{\infty}
\frac{(1-aq^{2n})}{(1-a)}\frac{(b,c,a/b;q)_n}{(aq/b,aq/c,bq;q)_n}
q^{\binom n2}\left(-\frac{aq}c\right)^n\\\times
(bq^{1+n},bq^{1-n}/a,c^2q^{2+n}/b,c^2q^{2-n}/ab;q^2)_{\infty}
\frac{(bcq/a,bq/c,bq^{1-n}/ad,bq^{1+n}/d;q)_N}
{(bq^{1-n}/a,bq^{1+n},bcq/ad,bq/cd;q)_N}\\
=\frac{(bcq/a,bq/c,bq/ad,bq/d;q)_N}
{(bq/a,bq,bcq/ad,bq/cd;q)_N}
\sum_{n=-\infty}^{\infty}
\frac{(1-aq^{2n})}{(1-a)}\frac{(b,c,ad/b,bq^{1+N}/d;q)_n}
{(aq/b,aq/c,bq/d,adq^{-N}/b;q)_n}\\\times
\frac{(aq^{-N}/b;q)_n}{(bq^{1+N};q)_n}
q^{\binom n2}\left(-\frac{aq}c\right)^n
(bq^{1+n},bq^{1-n}/a,c^2q^{2+n}/b,c^2q^{2-n}/ab;q^2)_{\infty}.
\end{multline*}
Equating the last expression with \eqref{ex14gl}, and performing the
substitution $d\mapsto bd/a$, we obtain the following bilateral quadratic
transformation:
\begin{Corollary}\label{cor4}
Let $a$, $b$, $c$, and $d$ be indeterminates, and let $N$ be a
nonnegative integer. Then
\begin{multline}\label{ex15gl}
\sum_{n=-\infty}^{\infty}
\frac{(1-aq^{2n})}{(1-a)}\frac{(b,c,d,aq^{1+N}/d,aq^{-N}/b;q)_n}
{(aq/b,aq/c,aq/d,dq^{-N},bq^{1+N};q)_n}
q^{\binom n2}\left(-\frac{aq}c\right)^n\\\times
(bq^{1+n},bq^{1-n}/a,c^2q^{2+n}/b,c^2q^{2-n}/ab;q^2)_{\infty}\\
=\frac{(q,q,aq,q/a,aq/bc,bq/c;q)_{\infty}}
{(bq,q/b,bq/a,aq/b,q/c,aq/c;q)_{\infty}}\,
(q,b^2q/a,c^2q^2/a,c^2q^2/b^2;q^2)_{\infty}\\\times
\frac{(bq/a,bq,cq/d,aq/cd;q)_N}{(bcq/a,bq/c,q/d,aq/d;q)_N}\,
{}_{10}\phi_9\!\left[\begin{matrix}bc/a,\,
q(bc/a)^{\frac 12},-q(bc/a)^{\frac 12},\\
(bc/a)^{\frac 12},-(bc/a)^{\frac 12},\end{matrix}\right.\\*
\left.\begin{matrix}bc/a,b/c,bd/a,c(q/a)^{\frac 12},-c(q/a)^{\frac 12},
bq^{1+N}/d,q^{-N}\\
q,c^2q/a,cq/d,b(q/a)^{\frac 12},-b(q/a)^{\frac 12},
cdq^{-N}/a,bcq^{1+N}/a\end{matrix}\,;q,
-\frac{cq}b\right].
\end{multline}
\end{Corollary}
Note that the bilateral sum on the left hand side of \eqref{ex15gl}
converges absolutely for $|q|<1$. See Theorem~\ref{qun} for a very similar
situation.

Now let us combine Theorem~\ref{qun} with Shukla's very-well-poised
${}_8\psi_8$ summation \eqref{shuklagl}.
We multiply both sides of \eqref{qungl} by
\begin{equation*}
\frac{(1-bcq^{2k}/a)}{(1-bc/a)}\frac{(c,c/a,d,e,fq,bcq/af;q)_k}
{(bq/a,bq,bcq/ad,bcq/ae,bc/af,f;q)_k}\left(\frac{b^2}{ade}\right)^k
\end{equation*}
and sum over all integers $k$. On the right hand side we obtain
\begin{multline}\label{ex16gl}
\frac{(q,q,aq,q/a,aq/bc,bq/c;q)_{\infty}}
{(bq,q/b,bq/a,aq/b,q/c,aq/c;q)_{\infty}}\,
(q,b^2q/a,c^2q^2/a,c^2q^2/b^2;q^2)_{\infty}\\\times
\sum_{k=-\infty}^{\infty}
\frac{(1-bcq^{2k}/a)}{(1-bc/a)}\frac{(bc/a,b/c,d,e,fq,bcq/af;q)_k}
{(q,c^2q/a,bcq/ad,bcq/ae,bc/af,f;q)_k}\\\times
\frac{(c^2q/a;q^2)_k}{(b^2q/a;q^2)_k}\left(-\frac{bc}{ade}\right)^k.
\end{multline}
On the left hand side we obtain
\begin{multline*}
\sum_{k=-\infty}^{\infty}
\frac{(1-bcq^{2k}/a)}{(1-bc/a)}\frac{(c,c/a,d,e,fq,bcq/af;q)_k}
{(bq/a,bq,bcq/ad,bcq/ae,bc/af,f;q)_k}\left(\frac{b^2}{ade}\right)^k\\\times
\sum_{n=-\infty}^{\infty}
\frac{(1-aq^{2n})}{(1-a)}\frac{(b,cq^k,aq^{-k}/b;q)_n}
{(aq/b,aq^{1-k}/c,bq^{1+k};q)_n}
q^{\binom n2}\left(-\frac{aq}c\right)^n\\\times
(bq^{1+n},bq^{1-n}/a,c^2q^{2+n}/b,c^2q^{2-n}/ab;q^2)_{\infty}\\
=\sum_{n=-\infty}^{\infty}
\frac{(1-aq^{2n})}{(1-a)}\frac{(b,c,a/b;q)_n}{(aq/b,aq/c,bq;q)_n}
q^{\binom n2}\left(-\frac{aq}c\right)^n\\\times
(bq^{1+n},bq^{1-n}/a,c^2q^{2+n}/b,c^2q^{2-n}/ab;q^2)_{\infty}\\\times
\sum_{k=-\infty}^{\infty}
\frac{(1-bcq^{2k}/a)}{(1-bc/a)}\frac{(cq^n,cq^{-n}/a,d,e,fq,bcq/af;q)_k}
{(bq^{1-n}/a,bq^{1+n},bcq/ad,bcq/ae,bc/af,f;q)_k}
\left(\frac{b^2}{ade}\right)^k.
\end{multline*}
Now the inner sum, provided $|b^2/ade|<1$,
can be evaluated by \eqref{shuklagl} and we obtain
\begin{multline*}
\sum_{n=-\infty}^{\infty}
\frac{(1-aq^{2n})}{(1-a)}\frac{(b,c,a/b;q)_n}{(aq/b,aq/c,bq;q)_n}
q^{\binom n2}\left(-\frac{aq}c\right)^n\\\times
(bq^{1+n},bq^{1-n}/a,c^2q^{2+n}/b,c^2q^{2-n}/ab;q^2)_{\infty}\\\times
\left(1-\frac{(1-c/b)(1-adq^n/b)(1-aeq^n/b)}
{(1-cq^n/f)(1-afq^n/b)(1-ade/b^2)}\right)
\frac{(1-fq^{-n}/c)(1-afq^n/b)}{(1-af/bc)(1-f)}\\\times
\frac{(q,bcq/a,aq/bc,bq/c,bq^{1-n}/ad,bq^{1-n}/ae,bq^{1+n}/d,bq^{1+n}/e,
bcq/ade;q)_{\infty}}{(bq^{1-n}/a,bq^{1+n},bcq/ad,bcq/ae,q^{1-n}/c,
aq^{1+n}/c,q/d,q/e,b^2q/ade;q)_{\infty}}\\
=\frac{(q,bcq/a,aq/bc,b/c,bq/ad,bq/ae,bq/d,bq/e,bcq/ade;q)_{\infty}}
{(bq/a,bq,bcq/ad,bcq/ae,q/c,aq/c,q/d,q/e,b^2/ade;q)_{\infty}}\\\times
\frac{bf}{ade(1-f)(1-af/bc)}
\sum_{n=-\infty}^{\infty}
\frac{(1-aq^{2n})}{(1-a)}\frac{(b,adq/b,aeq/b;q)_n}{(aq/b,bq/d,bq/e;q)_n}
q^{\binom n2}\left(-\frac b{de}\right)^n\\\times
(bq^{1+n},bq^{1-n}/a,c^2q^{2+n}/b,c^2q^{2-n}/ab;q^2)_{\infty}\\\times
\left(1-\frac{(1-cq^n/f)(1-afq^n/b)(1-ade/b^2)}
{(1-c/b)(1-adq^n/b)(1-aeq^n/b)}\right).
\end{multline*}
Equating the last expression with \eqref{ex16gl}, and performing the
simultaneous substitutions $c\mapsto e$, $d\mapsto bc/a$ and $e\mapsto bd/a$, 
we obtain the following bilateral quadratic transformation:
\begin{Corollary}\label{cor5}
Let $a$, $b$, $c$, $d$, $e$, and $f$ be indeterminates. Then
\begin{multline}\label{ex17gl}
\sum_{n=-\infty}^{\infty}
\frac{(1-aq^{2n})}{(1-a)}\frac{(b,cq,dq;q)_n}{(aq/b,aq/c,aq/d;q)_n}
q^{\binom n2}\left(-\frac {a^2}{bcd}\right)^n\\\times
(bq^{1+n},bq^{1-n}/a,e^2q^{2+n}/b,e^2q^{2-n}/ab;q^2)_{\infty}\\\times
\left(1-\frac{(1-eq^n/f)(1-afq^n/b)(1-cd/a)}
{(1-e/b)(1-cq^n)(1-dq^n)}\right)\\
=\frac{bcd(1-f)(1-af/be)\,(q,aq,q/a,eq/c,eq/d,aq/bc,aq/bd,a/cd;q)_{\infty}}
{af(1-b/e)\,(q/b,aq/b,beq/a,q/c,q/d,aq/c,aq/d,aeq/bcd;q)_{\infty}}\\\times
(q,b^2q/a,e^2q^2/a,e^2q^2/b^2;q^2)_{\infty}\,
{}_{10}\phi_9\!\left[\begin{matrix}be/a,\,
q(be/a)^{\frac 12},-q(be/a)^{\frac 12},\\
(be/a)^{\frac 12},-(be/a)^{\frac 12},\end{matrix}\right.\\*
\left.\begin{matrix}b/e,bc/a,bd/a,fq,beq/af,
e(q/a)^{\frac 12},-e(q/a)^{\frac 12}\\
e^2q/a,eq/c,eq/d,be/af,f,b(q/a)^{\frac 12},-b(q/a)^{\frac 12}
\end{matrix}\,;q,
-\frac{ae}{bcd}\right],
\end{multline}
provided $|ae/bcd|<1$.
\end{Corollary}
Note that we have applied analytic continuation to obtain the convergence
condition in Corollary~\ref{cor5}.

If we let multiply both sides of \eqref{ex17gl} by $f$ and then let $f\to 0$,
we obtain
\begin{multline}\label{ex18gl}
\sum_{n=-\infty}^{\infty}
\frac{(1-aq^{2n})}{(1-a)}\frac{(b,c,d;q)_n}{(aq/b,aq/c,aq/d;q)_n}\,
q^{\binom n2}\left(-\frac {a^2q}{bcd}\right)^n\\\times
(bq^{1+n},bq^{1-n}/a,e^2q^{2+n}/b,e^2q^{2-n}/ab;q^2)_{\infty}\\
=\frac{(q,aq,q/a,eq/c,eq/d,aq/bc,aq/bd,aq/cd;q)_{\infty}}
{(q/b,aq/b,beq/a,q/c,q/d,aq/c,aq/d,aeq/bcd;q)_{\infty}}\\\times
(q,b^2q/a,e^2q^2/a,e^2q^2/b^2;q^2)_{\infty}\,
{}_8\phi_7\!\left[\begin{matrix}be/a,\,
q(be/a)^{\frac 12},-q(be/a)^{\frac 12},\\
(be/a)^{\frac 12},-(be/a)^{\frac 12},\end{matrix}\right.\\
\left.\begin{matrix}b/e,bc/a,bd/a,
e(q/a)^{\frac 12},-e(q/a)^{\frac 12}\\
e^2q/a,eq/c,eq/d,b(q/a)^{\frac 12},-b(q/a)^{\frac 12}
\end{matrix}\,;q,
-\frac{aeq}{bcd}\right],
\end{multline}
provided $|aeq/bcd|<1$. Now if we specialize \eqref{ex18gl}
by letting $d\to a/c$ and $e\to c$, the ${}_8\phi_7$ series
on the right hand side can be evaluated by the $q$-Watson summation
in \eqref{qwatsongl} and we obtain the following bilateral
quadratic summation formula:
\begin{Corollary}\label{cor6}
Let $a$, $b$, and $c$  be indeterminates. Then
\begin{multline}\label{ex19gl}
\sum_{n=-\infty}^{\infty}
\frac{(1-aq^{2n})}{(1-a)}\frac{(b,c,a/c;q)_n}{(aq/b,aq/c,cq;q)_n}
q^{\binom n2}\left(-\frac {aq}{b}\right)^n\\\times
(bq^{1+n},bq^{1-n}/a,c^2q^{2+n}/b,c^2q^{2-n}/ab;q^2)_{\infty}\\
=\frac{(q,q,aq,q/a,aq/bc,cq/b;q)_{\infty}}
{(q/b,aq/b,cq,q/c,aq/c,cq/a;q)_{\infty}}\,
(bq/c,bcq/a,cq^2/b,c^3q^2/b;q^2)_{\infty},
\end{multline}
provided $|cq/b|<1$.
\end{Corollary}

\subsection{A consequence of Theorem~\ref{shn}}
For illustration, we combine Theorem~\ref{shn} with
Shukla's very-well-poised ${}_8\psi_8$ summation.
We multiply both sides of \eqref{ex4gl} (which is equivalent to
\eqref{shngl}) by
\begin{equation*}
\frac{(1-aq^{2k})}{(1-a)}\frac{(a/c,b,f,g,v,aq^2/v;q)_k}
{(cq,aq/b,aq/g,aq/g,aq/v,v/q;q)_k}\left(\frac{ac}{bfg}\right)^k
\end{equation*}
and sum over all integers $k$. On the right hand side we obtain
\begin{multline}\label{ex20gl}
\frac {(1-u/q)(1-aq/u)(q,aq/d,aq/e,q/d,q/e,cq/b,aq/bc,bcq/a,ab/cde;q)_{\infty}}
{(b-c)(a/cd,a/ce,bq/d,bq/e,q/b,aq/b,cq,cq/a,aq/de;q)_{\infty}}\\\times
{}_{10}\psi_{10}\!\left[\begin{matrix}
q\sqrt{a},-q\sqrt{a},d,e,f,g,u,aq^2/u,v,aq^2/v\\
\sqrt{a},-\sqrt{a},aq/d,aq/e,aq/f,aq/g,aq/u,u/q,aq/v,v/q\end{matrix}\,;q,
\frac{a^2}{defgq}\right].
\end{multline}
On the left hand side we obtain
\begin{multline*}
\sum_{k=-\infty}^{\infty}
\frac{(1-aq^{2k})}{(1-a)}\frac{(a/c,b,f,g,v,aq^2/v;q)_k}
{(cq,aq/b,aq/g,aq/g,aq/v,v/q;q)_k}\left(\frac{ac}{bfg}\right)^k\\\times
\sum_{n=-\infty}^{\infty}\frac{(1-bcq^{2n}/a)}{(1-bc/a)}
\frac{(bq^k,bq^{-k}/a,cdq/a,ceq/a;q)_n}
{(cq^{1-k}/a,cq^{1+k},bq/d,bq/e;q)_n}
\left(\frac a{de}\right)^n\\\times
\left(1-\frac{(1-cde/ab)(1-cq^{n+1}/u)(1-cuq^{n-1}/a)}
{(1-c/b)(1-cdq^n/a)(1-ceq^n/a)}\right)\\
=\sum_{n=-\infty}^{\infty}
\frac{(1-bcq^{2n}/a)}{(1-bc/a)}
\frac{(b,b/a,cdq/a,ceq/a;q)_n}{(cq/a,cq,bq/d,bq/e;q)_n}
\left(\frac a{de}\right)^n\\\times
\left(1-\frac{(1-cde/ab)(1-cq^{n+1}/u)(1-cuq^{n-1}/a)}
{(1-c/b)(1-cdq^n/a)(1-ceq^n/a)}\right)\\\times
\sum_{k=-\infty}^{\infty}
\frac{(1-aq^{2k})}{(1-a)}\frac{(aq^{-n}/c,bq^n,f,g,v,aq^2/v;q)_k}
{(cq^{1+n},aq^{1-n}/b,aq/g,aq/g,aq/v,v/q;q)_k}
\left(\frac{ac}{bfg}\right)^k.
\end{multline*}
Now the inner sum, provided $|ac/bfg|<1$, can be evaluated by Shukla's
very-well-poised ${}_8\psi_8$ summation \eqref{shuklagl} and we obtain
\begin{multline*}
\sum_{n=-\infty}^{\infty}
\frac{(1-bcq^{2n}/a)}{(1-bc/a)}
\frac{(b,b/a,cdq/a,ceq/a;q)_n}{(cq/a,cq,bq/d,bq/e;q)_n}
\left(\frac a{de}\right)^n\\\times
\left(1-\frac{(1-cde/ab)(1-cq^{n+1}/u)(1-cuq^{n-1}/a)}
{(1-c/b)(1-cdq^n/a)(1-ceq^n/a)}\right)\\\times
\left(1-\frac{(1-b/c)(1-bfq^n/a)(1-bgq^n/a)}
{(1-bq^{1+n}/v)(1-bvq^{n-1}/a)(1-bfg/ac)}\right)
\frac{(1-vq^{-1-n}/b)(1-bvq^{n-1}/a)}{(1-v/aq)(1-v/q)}\\\times
\frac{(q,aq,q/a,cq/b,aq^{1-n}/bf,aq^{1-n}/bg,cq^{1+n}/f,cq^{1+n}/g,
aq/fg;q)_{\infty}}
{(aq^{1-n}/b,cq^{1+n},aq/f,aq/g,q^{1-n}/b,cq^{1+n}/a,
q/f,q/g,acq/bfg;q)_{\infty}}\\
=\frac{bv(q,aq,q/a,cq/b,aq/bf,aq/bg,cq/f,cq/g,aq/fg;q)_{\infty}}
{aq(1-v/aq)(1-v/q)
(aq/b,cq,aq/f,aq/g,q/b,cq/a,q/f,q/g,acq/bfg;q)_{\infty}}\\\times
\sum_{n=-\infty}^{\infty}
\frac{(1-bcq^{2n}/a)}{(1-bc/a)}
\frac{(cdq/a,ceq/a,bfq/a,bgq/a;q)_n}{(bq/d,bq/e,cq/f,cq/g;q)_n}
\left(\frac a{deq}\right)^n\\\times
\left(1-\frac{(1-cde/ab)(1-cq^{n+1}/u)(1-cuq^{n-1}/a)}
{(1-c/b)(1-cdq^n/a)(1-ceq^n/a)}\right)\\\times
\left(1-\frac{(1-bq^{1+n}/v)(1-bvq^{n-1}/a)(1-bfg/ac)}
{(1-b/c)(1-bfq^n/a)(1-bgq^n/a)}\right).
\end{multline*}
Equating the last expression with \eqref{ex20gl},
we obtain the following transformation for a particular
very-well-poised ${}_{10}\psi_{10}$ series:
\begin{Corollary}\label{cor7}
Let $a$, $b$, $c$, $d$, $e$, $f$, $g$, $u$, and $v$ be indeterminates. Then
\begin{multline}\label{ex21gl}
{}_{10}\psi_{10}\!\left[\begin{matrix}
q\sqrt{a},-q\sqrt{a},d,e,f,g,u,aq^2/u,v,aq^2/v\\
\sqrt{a},-\sqrt{a},aq/d,aq/e,aq/f,aq/g,aq/u,u/q,aq/v,v/q\end{matrix}\,;q,
\frac{a^2}{defgq}\right]\\
=\frac{(b-c)(c-b)}
{(1-u/q)(1-aq/u)(1-v/q)(1-aq/v)}\\\times
\frac{(aq,q/a,a/cd,a/ce,aq/de,aq/bf,aq/bg,aq/fg,
bq/d,bq/e,cq/f,cq/g;q)_{\infty}}
{(q/d,aq/d,q/e,aq/e,q/f,aq/f,q/g,aq/g,
aq/bc,bcq/a,ab/cde,acq/bfg;q)_{\infty}}\\\times
\sum_{n=-\infty}^{\infty}
\frac{(1-bcq^{2n}/a)}{(1-bc/a)}
\frac{(cdq/a,ceq/a,bfq/a,bgq/a;q)_n}{(bq/d,bq/e,cq/f,cq/g;q)_n}
\left(\frac a{deq}\right)^n\\\times
\left(1-\frac{(1-cde/ab)(1-cq^{n+1}/u)(1-cuq^{n-1}/a)}
{(1-c/b)(1-cdq^n/a)(1-ceq^n/a)}\right)\\\times
\left(1-\frac{(1-bq^{1+n}/v)(1-bvq^{n-1}/a)(1-bfg/ac)}
{(1-b/c)(1-bfq^n/a)(1-bgq^n/a)}\right),
\end{multline}
provided $|a^2/defgq|<1$ and $|a/deq|<1$.
\end{Corollary}
As usual, we have applied analytic continuation to establish
the conditions of convergence.

We could also combine Theorem~\ref{shn} with \eqref{milnekm2gl}
or with other identities. With Corollary~\ref{cor7}, we just wanted
to provide one example from the many possibilities of deducing
a transformation from Theorem~\ref{shn}.

\end{document}